\documentclass{egpubl}
\usepackage{eg2025}
 
\ConferenceSubmission   % uncomment for Conference submission
\usepackage[T1]{fontenc}
\usepackage{dfadobe}  

\usepackage{cite}  % comment out for biblatex with backend=biber
\BibtexOrBiblatex
\electronicVersion
\PrintedOrElectronic
\ifpdf \usepackage[pdftex]{graphicx} \pdfcompresslevel=9
\else \usepackage[dvips]{graphicx} \fi

\usepackage{egweblnk} 

\title{Automatic Ply Partitioning for Laminar Composite Process Planning}

\author[E. Garner \& A.\, M. Mirzendehdel]
{\parbox{\textwidth}{\centering E. Garner$^{1,2}$\orcid{0009-0002-2323-5833}
        and A.\, M. Mirzendehdel$^{2}$\orcid{0000-0002-4407-1877} 
        }
        \\
{\parbox{\textwidth}{\centering $^1$Université de Lorraine, CNRS, Inria, LORIA, Nancy, France\\
         $^2$SRI International, Palo Alto, California
       }
}
}

\usepackage{enumitem}
\usepackage{amsmath}
\usepackage{amssymb}
\usepackage{algorithm2e}
\newcommand{\sign}{\mathsf{sign}}
\newcommand{\minimize}{\mathop{ \text{minimize}}} 

\begin{document}

\teaser{
 \includegraphics[width=1\linewidth]{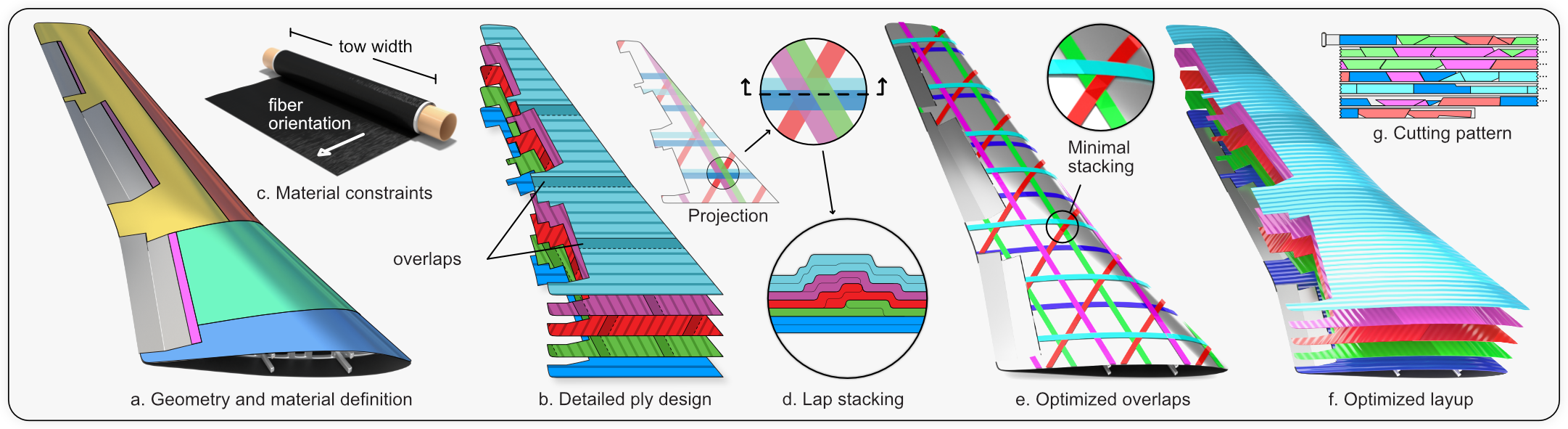}
 \centering
  \caption{Laminar composite ply partitioning pipeline: Geometry and material definition specify part geometry and zonal mechanical properties (a). Software-generated detailed ply designs embody design intent (b), but cannot account for constraints on material dimensions (c). For large-scale parts, the final step of the detailed ply design stage involves manually partitioning plies into manufacturable sub-plies. This step introduces overlaps across adjacent sub-plies, potentially leading to defects such as overlap stacking (d). The proposed algorithm optimizes this partitioning step to minimize adverse effects on part quality (e) and optimize cutting patterns (g).}
\label{fig:teaser}
}

\maketitle
%-------------------------------------------------------------------------
\begin{abstract}

This work introduces an automated ply partitioning strategy for large-scale laminar composite manufacturing. It specifically targets the problem of fabricating large plies from available spooled materials, while minimizing the adverse effects on part quality. The proposed method inserts fiber-aligned seams sequentially until each resulting sub-ply can be manufactured from available materials, while simultaneously enforcing constraints to avoid quality issues induced by the stacking of seams across multiple plies. Leveraging the developable nature of individual plies, the partitioning problem is cast as a sequence of one-dimensional piecewise linear optimization problems, thus allowing for efficient local optimization via linear programming. We experimentally demonstrate that coupling the local search with a greedy global search produces the same results as an exhaustive search. The resulting automated method provides an efficient and robust alternative to the existing trial-and-error approach, and can be readily integrated into state-of-the-art composite design workflows. In addition, this formulation enables the inclusion of common constraints regarding laminate thickness tolerance, sub-ply geometry, stay-out zones, material wastage, etc. The efficacy of the proposed method is demonstrated through its application to the surface of an airplane wing and to the body panels of an armored vehicle, each subject to various performance and manufacturing-related geometric constraints.

\end{abstract}  
%-------------------------------------------------------------------------
\section{Introduction}
\label{sec:intro}
Fiber-reinforced laminar composites are increasingly popular in the design of large-scale structures such as wind turbine blades, aircraft fuselages, and other industrial equipment where their excellent strength to weight ratio, design flexibility, and environmental stability offer compelling advantages over traditional materials. However, the complex design and manufacturing pipeline associated with laminar composites remains an impediment to wider adoption~\cite{Grand_View_Research_2024}. 

\begin{figure}
    \centering
    \includegraphics[width=\linewidth]{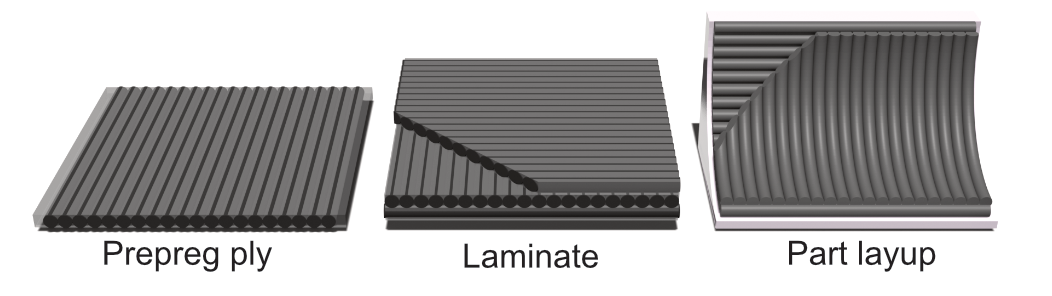}
    \caption{Laminar composite elements: Prepreg ply, composed for matrix-embedded unidirectional fibers is layered in specific orientations during layup to form a laminate with a desired shape and mechanical properties.}
    \label{fig:laminate_illustation}
\end{figure}

In broad terms, laminar composite manufacturing involves layering oriented fibers, impregnating them with resin, and consolidating the assembly into a unified part. Although a wide variety of manufacturing methods exist, large-scale components are most commonly fabricated following the \emph{prepreg} method. This method streamlines the process by utilizing spools of axially-aligned reinforcing fibers pre-impregnated with semi-cured resin. Plies are cut from the prepreg material, arranged by hand into a layup, and cured to create the final consolidated laminate part~\cite{staab2015laminar} (Fig~\ref{fig:laminate_illustation}).

\begin{figure}
    \centering
    \includegraphics[width=\linewidth]{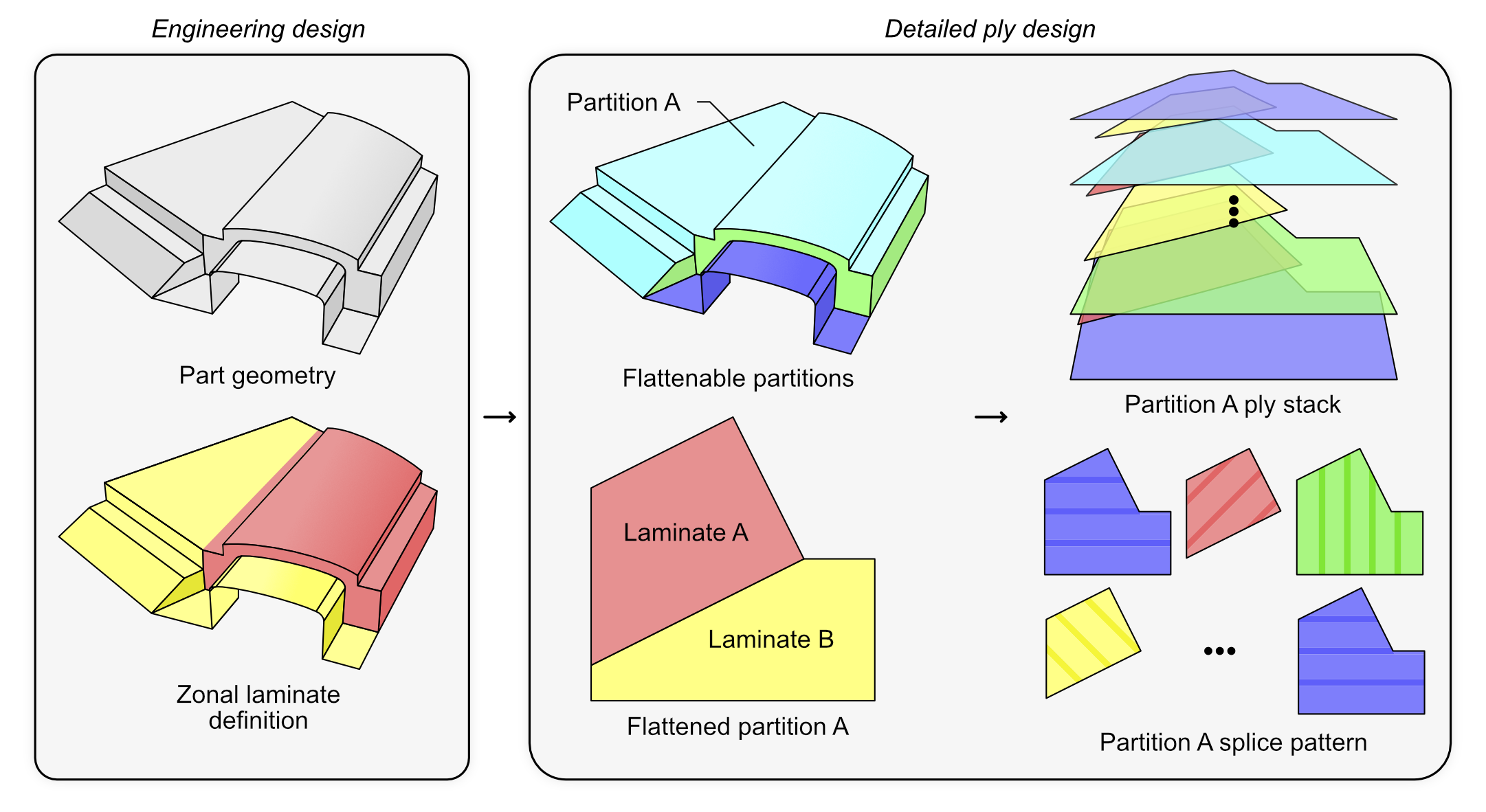}
    \caption{Design pipeline: Part geometry and composite properties are defined in the engineering design stage. The part is partitioned into flattenable regions and detailed ply designs are generated for each region in the detailed ply design stage.}
    \label{fig:design_pipeline}
\end{figure}

Prior to manufacturing, a comprehensive process planning pipeline translates the intended design into a set of cutting patterns, tooling requirements, process parameters, and assembly instructions. The pipeline can be broadly divided into three stages: \emph{engineering design}, \emph{detailed ply design}, and \emph{manufacturing design}~\cite{fibersim_siemens}. In the first stage, the part geometry, spatially-varying mechanical properties, and specific materials are defined. In the second stage, the geometry is partitioned into developable sections, and detailed ply layups are generated to achieve the desired properties, while incorporating engineering specifications. These two stages are illustrated in Fig.~\ref{fig:design_pipeline}. In the third stage, the required tooling is designed and prepreg cutting patterns, layup sequences, and process parameters are generated.

The inherent complexity of the process planning pipeline has prompted the development of specialized CAD systems such as NX FiberSIM~\cite{fibersim_siemens}, and CATIA Composite Design~\cite{noauthor_catia_2023} to streamline and automated the process. Yet, certain tasks continue to rely heavily on labor-intensive trial-and-error methods. One such task is the partitioning of plies to account for prepreg spool dimensions.

\begin{figure}
    \centering
    \includegraphics[width=\linewidth]{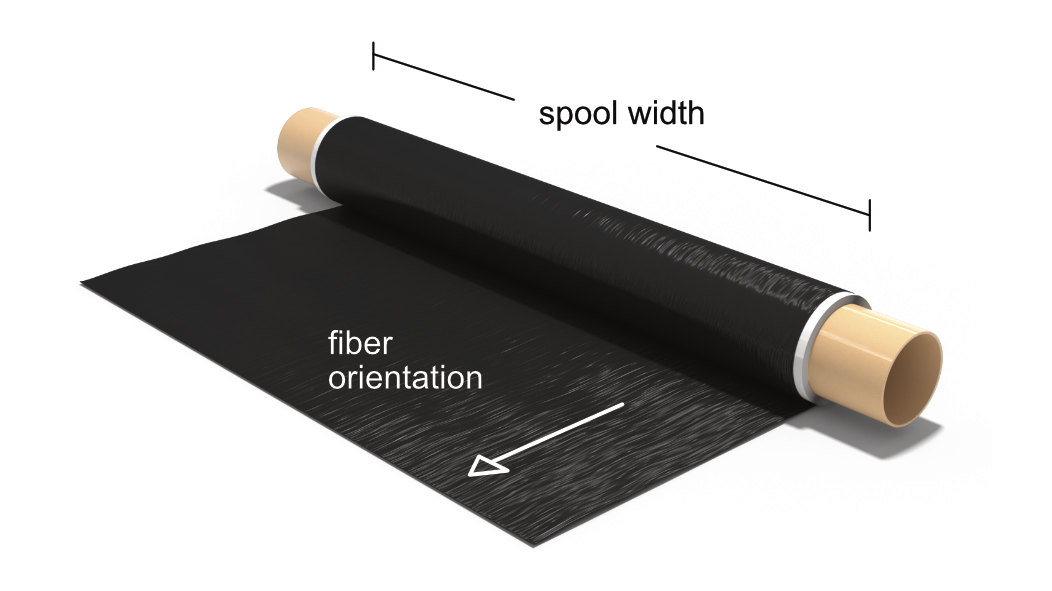}
    \caption{Prepreg spool illustration.}
    \label{fig:spool}
\end{figure}

Prepreg spools, from which plies are cut, are wound along the fiber direction, effectively limiting the dimensions of each ply in the fiber-transverse direction (Fig.~\ref{fig:spool}). For large-scale components, it is often necessary to partition plies into smaller sub-plies such that they can be fabricated from available materials. The resulting partition seams between adjacent sub-plies introduce potential weak spots, which can be exacerbated by positioning inaccuracy. Manufacturers mitigate this risk by introducing slight overlaps between adjacent sub-plies. However, these overlaps introduce yet another deviation from the original design intent, and can have a compounding effect when they intersect one another through the thickness of the laminate (Fig.~\ref{fig:teaser}d); thus their placement requires careful consideration \cite{mehdikhani2019voids}. Existing tools can recommend partitioning when plies exceed material dimensions and offer simple seam staggering options, but to the best of our knowledge, offer little in terms of design automation, or optimization.

The present work targets this specific task, and introduces the first automated ply partitioning strategy to meet spool width constraints, while minimizing the adverse effects on part quality. Our core contribution is a ply partitioning algorithm that ensures manufacturability from available materials, while incorporating common design considerations such as spatially-varying geometric tolerances, ease of assembly, and material wastage. 

The geometry and material definition, as well as the prior steps of the detailed ply design stage fall outside of the scope of this paper, and are treated as inputs to the algorithm.

\section{Related Work}
\label{sec:related_work}

The ply partitioning problem is one of dividing a set of continuous domains (plies) into discrete subdomains (sub-plies), subject to geometric constraints. Although the particular design constraints are application-specific, geometry partitioning and parametrization are well-studied topics, with many applications in computer-aided design. Here, we present a brief overview of important methods and applications in which process planning has been optimized to suit manufacturing constraints.

Ply partitioning can be thought of as performing cuts on surfaces in 3D. This is reminiscent of subtractive manufacturing methods such as water jet, laser, or wire cutting, in which chunks of material are sequentially removed from a workpiece until the desired shape is revealed. In this context, a general surface is approximated as a piecewise ruled surfaces resulting from straight cuts. A large body of research focuses on fitting general surfaces to well-behaved ones such as piecewise ruled or developable surfaces amenable to specific manufacturing methods\cite{pottmann1999approximation, elber19975, subag2006piecewise}. In the related field of computer numerically control (CNC) machining, process planning involves deconstructing the workpiece into a series of cutting paths, each of which removes a controlled amount of material. Here, partitioning methods must take into account not only the workpiece geometry, but also the machine and cutting tool specifications~\cite{elber1993tool, kim2015precise}. Similarly, in additive manufacturing, digital designs are first partitioned into discrete slabs, and then into tool paths. In this context, specialized methods have been proposed to account for constraints on geometric tolerances~\cite{masood2000part, smith2002optimal, alexa2017optimal}, mechanical performance~\cite{ren2019thermo, tura2022characterization}, surface quality~\cite{thrimurthulu2004optimum, canellidis2006pre, delfs2016optimized, jin2017optimization}, build time~\cite{thrimurthulu2004optimum, canellidis2006pre, delfs2016optimized}, cost~\cite{alexander1998part}, and build volume~\cite{luo2012chopper}.

More closely related to our application, digital methods 
for garment design include strategies for decomposing 3D surfaces into flattenable patches~\cite{wang2002surface, poranne2017autocuts}, sewing pattern optimization for different body shapes~\cite{wolff2023designing, meng2012flexible}, seam durability~\cite{montes2020computational}, and aesthetics~\cite{kwok2015styling}. These methods build on techniques for optimal surface unwrapping via energy minimization~\cite{poranne2017autocuts}, edge clustering for aesthetic purposes~\cite{bruno2009you}, and packing algorithms for waste minimization~\cite{nee1986designing}.
Recently, Pietroni \emph{et al.} proposed an interactive tool for creating sewing patterns from 3D garment models~\cite{pietroni2022computational}. Their work aimed to optimize the position of seams to balance material, aesthetic, and manufacturability-related constraints. Similar to process planning for laminar composites, they partition 3D surfaces into flattenable patches that optimize the trade-off between distortion and total seam length. Seam placement takes into account both the 3D configuration for aesthetic considerations, as well as the 2D configuration for manufacturing constraints.

In the context of laminar composite design, several important distinctions from garment design apply: \emph{1. prepreg is extremely rigid}, and \emph{2. seams are aligned with fiber orientation}. \emph{1} ensures that parts are designed from the outset to avoid non-developable surfaces. Thus, a solution can be developed in the 2D domain. \emph{2} means that seams can be defined by a single variable corresponding to an offset from the origin.

These attributes allow for efficient linear optimization strategies. Specifically, the proposed algorithm casts the problem as a sequence of seam insertion problems defined as piecewise linear systems, and solves it through a combination of linear programming and heuristic tree search.
\section{Overview}
\label{sec:overview}

The proposed algorithm takes as input a set of ply designs, together with user-defined manufacturing constraints. The output is a set of sub-plies, as well as a cutting pattern, optimized according to the provided manufacturing constraints.

\subsection{Manufacturing constraints}
The key manufacturing constraint, and the reason partitioning is necessary in the first place, relates to the width of the prepreg spool. Since spools are wound in the fiber direction as shown in Fig.~\ref{fig:spool}, the spool width constrains the size of each ply in the fiber-transverse direction. The simplest partitioning strategy therefore consists of splitting the flattened representation of each ply along straight lines in the fiber direction. In theory, seams need not be fiber-aligned, nor straight. However, standard practice in industry is to place seams predominantly in the fiber direction. This minimizes the induced strength drop down and maximizes surface coverage for a single sub-ply.

Additional manufacturing constraints restrict seam placement to minimize adverse effects on as-manufactured part quality. These constraints are described below.

\subsubsection{Overlaps}
As mentioned above, manufacturers introduce small overlaps between sub-plies to avoid mechanical strength drop-down. This practice, however, potentially introduces a new issue, since overlaps are locally double thickness. While the effect of a single overlap is generally negligible, issues do arise when overlaps across multiple stacked plies intersect. In this instance, local thickness tolerances and overall surface quality may be significantly affected. It is therefore critical to avoid overlap stacking in the partitioning stage.

\subsubsection{Stay-out zones}
Another common constraint prohibits placement of seams in specific areas, referred to as stay-out zones, where mechanical performance or thickness tolerance is critical. These zones effectively break up the originally continuous design space for each seam offset by restricting ranges which would result in stay-out intersection.

\subsubsection{Prepreg cutting patterns}
In addition to constraints on seam placement, designers are concerned with the shape of the resulting sub-plies. Specifically, small or flimsy sub-plies should be avoided in order to reduce the risk of inaccurate placement during hand layup, transport, or part consolidation. Moreover, small sub-plies, and specifically narrow sub-plies, result in more sparse prepreg cutting patterns, and therefore greater material wastage. Such constraints can be thought of in terms of minimum offsets between seams and the flattened part boundary.

\subsubsection{Production cost and design fidelity}
Prepreg spools are offered in a range of widths, with the cost per unit area generally increasing as the width increases. On the one hand, using narrower prepreg reduces material costs; while on the other, wider prepreg leads to less partitioning and better design fidelity. Designers must carefully weigh the trade-off between production cost and design fidelity. 

\subsection{Manual partitioning}
\begin{figure}
    \centering
    \includegraphics[width=1.0\linewidth]{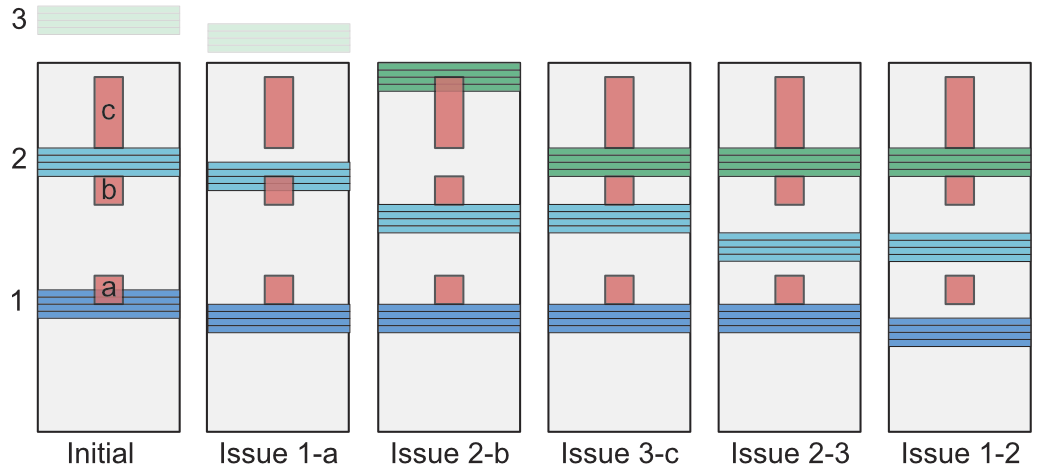}
    \caption{Manual partitioning operation sequence beginning from a spool width-driven pattern, and correcting issues sequentially until a viable design is found.}
    \label{fig:manual}
\end{figure}

Manual partitioning, though often trivial, frequently requires several time-consuming iterations. Consider the simple single-orientation example in Fig.~\ref{fig:manual}. A reasonable first attempt involves staggering the seams on each ply by the overlap width and placing each seam group as far apart as possible, according to the spool width. In this case, starting from the bottom, this results in seam group $1$ ends up intersecting stay-out zone $a$. In correcting this issue, seam group $2$ ends up intersecting stay-out zone $b$ due to the spool width constraint. Correcting this, in turn, requires the introduction of seam group $3$, which now intersects stay-out zone $c$. Correcting this causes seam groups $2$ and $3$ to be too close to one another, violating the minimum sub-ply width constraint. Two more iterations are required before a viable design is found. Naturally, this process can become very tedious with more complex examples, and there is no guarantee that it will produce a viable solution, whether or not one exists.

\subsection{Automated solution strategy}
Assuming plies are partitioned sequentially, and seams within each ply are placed to the right of the previous, then the insertion of a new seam presents a local optimization problem. If there exists a solution, the seam can be inserted, and the procedure can continue placing new seams until the next seam falls outside of the ply polygon, at which point, the process is repeated for the next ply.
Once all plies have been partitioned in this way, a viable solution has been found and the procedure terminates. Intuitively, a greedy strategy such as this may fail to find a viable solution; and a more exhaustive search strategy seems warranted \cite{cormen2022introduction}. Here, if the optimization objective is defined as maximizing the total distance covered by the sub-plies produced from a given set of seams, we observe that the problem is piecewise linear. A reasonable strategy is to search each linear subspace, exploring the design tree heuristically through a beam search. This solution offers a compromise between the efficiency of the greedy search and the robustness of an exhaustive one. We will demonstrate empirically, however, that the beam search does not perform better than the greedy search, and therefore seems unnecessary.

\section{Method}
\label{sec:method}

\begin{figure}
    \centering
    \includegraphics[width=\linewidth]{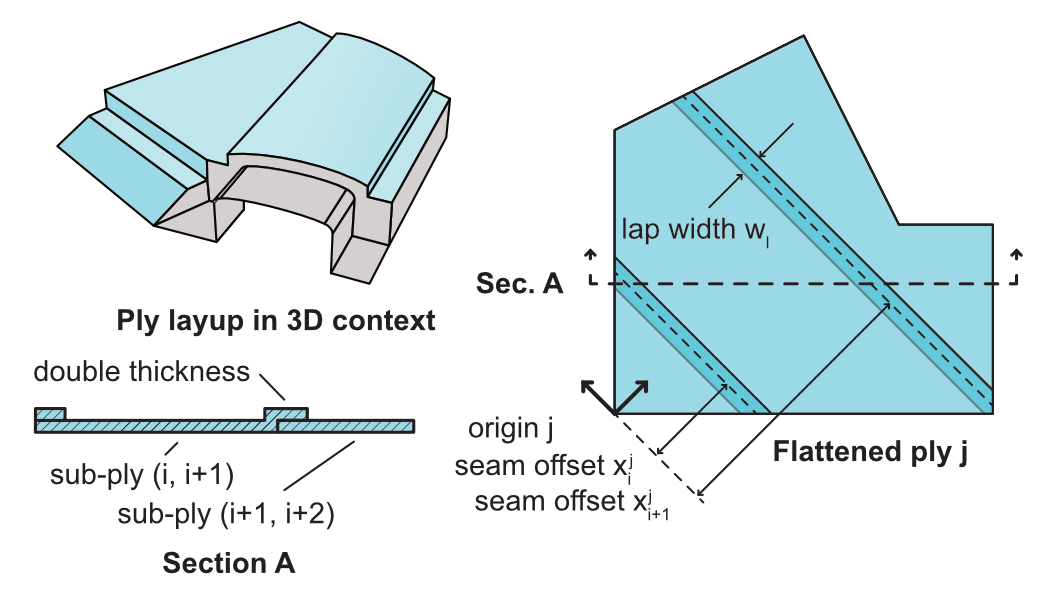}
    \caption{Example partitioned ply and key terminology. The 2D flattened ply representation is divided along straight seams in the fiber direction. Overlaps extend each resulting sub-ply symmetrically about the seam. Seams are defined in terms of their offset from a ply-specific origin.}
    \label{fig:terminology}
\end{figure}

The proposed design strategy takes as input a stack of flattened ply polygons, along with the spool width $w_s$ and other user-defined manufacturing constraints.
Seams are defined in the 2D space by line equations ${\mathbf{r} = \mathbf{r}_0 + \mathbf{d}_i t}$, where $\mathbf{d}_i$ is the fiber direction of the corresponding ply $i$. We define the design variable associated with seam $i$ of ply $j$ as its offset $x_i^j$ from the ply origin point $O_j$, as shown in Fig.~\ref{fig:terminology}. For each ply, the origin $O_j$ is chosen to ensure ensure that all offsets remain positive-valued.

Next, we represent the width of a sub-ply as the distance between adjacent seams on a given ply as ${x_{i+1}^j-x_{i}^j + w_l}$, where the overlap width $w_l$ term accounts for the fact that sub-plies are extended on either side by $\frac{w_l}{2}$ to create the overlap. The spool width constraint can then be expressed as:
\begin{equation}
    x_{i+1}^j-x_{i}^j+w_l \leq w_{s}.
    \label{eq:max_ply_width}
\end{equation}

\subsection{Overlap stacking}
\begin{figure}
    \centering
    \includegraphics[width=1.0\linewidth]{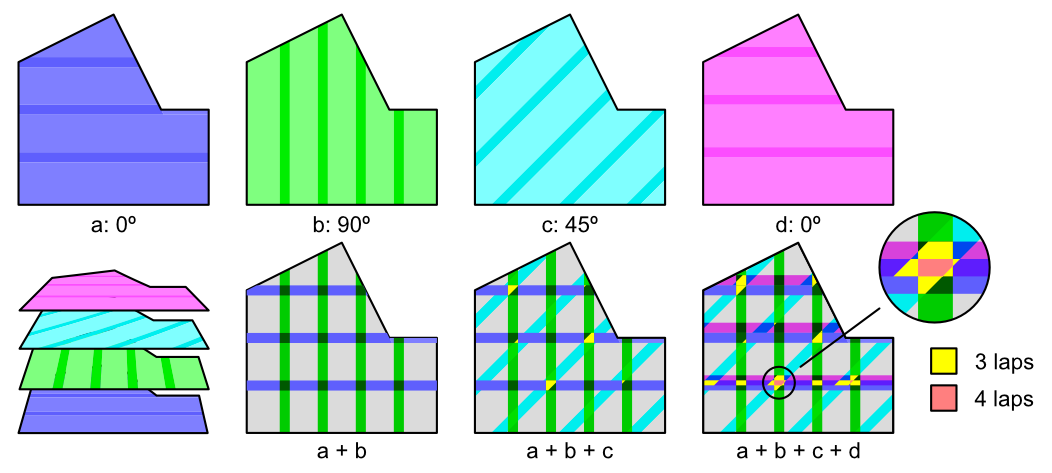}
    \caption{Top: Example partition pattern for a four-ply layup. Bottom: Overlap projection highlighting the compounding effect of intersecting seams.}
    \label{fig:overlap_illustration}
\end{figure}

When two or more seams intersect across multiple stacked plies as shown in Fig.~\ref{fig:overlap_illustration},
the thickness increase along the overlaps is compounded, producing an uneven surface. Intersections between two overlaps are largely unavoidable when the dimensions of the flattened plies are much greater than the available spool width. However, their effect is generally small over the full thickness of the laminate. Consider, for example, a laminate made up of 50 plies. A two-overlap stack would result in a thickness variation of 4\%; but as the number of intersecting overlaps increases, the thickness variation approaches 100\%. In practice, manufacturers often define a maximum tolerance for thickness variation $\delta$:
\begin{equation}
    \delta = \frac{t_{max} - t_{min}}{t_{min}},
\end{equation}
which can equivalently be written as the ratio of the maximum number of intersecting overlaps $N$ over the total number of plies $M$.
If one were to prohibit intersections of $N$ overlaps using only linear constraints, then $N$ constraints would be required.
Instead, constraints can be defined in terms of a smaller number of overlaps $n$ over a set of ply bundles, each made up of $m$ plies. In the worst case, intersections from each ply bundle may stack up, resulting in the original thickness variation of $N/M$. With this approach, the required number of plies in each bundle is:
\begin{equation}
    m = \Bigg{\lceil} \frac{Mn}{N} \Bigg{\rceil}.
\end{equation}
With two-overlap stacks being unavoidable, the smallest value of $n$ is 2, and constraints must be defined to avoid intersections between all possible combinations of $n + 1 = 3$ overlaps. Two cases may arise: either at least two seams are parallel, or all seams are non-parallel. In the first case, the intersecting area is non-zero when the parallel overlaps intersect, \emph{i.e.}, when the following inequality is violated:

\begin{equation}
\label{eq:lap_stacking_parallel}
    |x_j^* - x_i^*| \geq w_l.
\end{equation}
where $i$ and $j$ are parallel seams, and $^*$ specifies that they may be from any ply.

\begin{figure}
    \centering
    \includegraphics[width=1.0\linewidth]{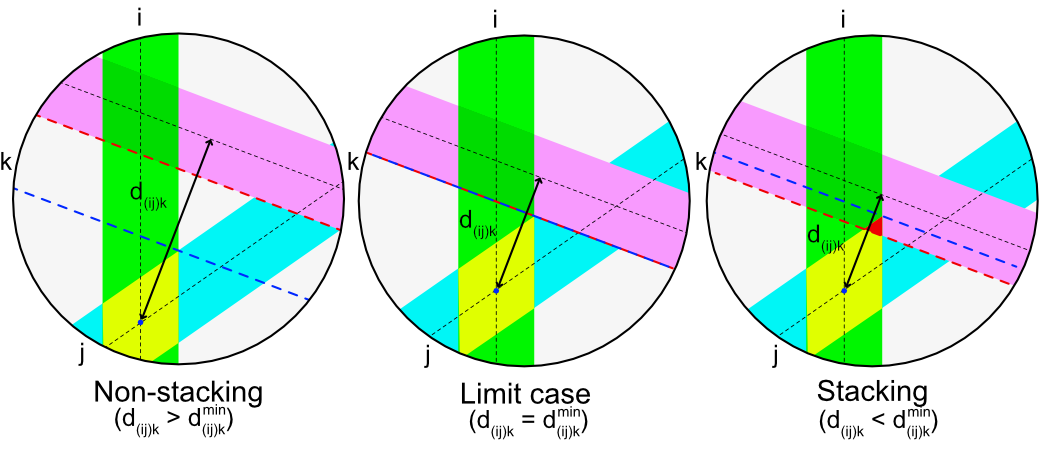}
    \caption{Triple overlap stacking condition between non-parallel overlaps.}
    \label{fig:triple_overlap_illustration}
\end{figure}

In the case of three non-parallel seams, three constraints are required. Consider the situation in Fig.~\ref{fig:triple_overlap_illustration} in which the triple intersection area is highlighted in red. Triple intersection occurs when the double intersection between $i$ and $j$, shown in yellow, intersects $k$. To avoid intersection, the distance between the projection of the $ij$ intersecting area in the direction normal to $k$ must be greater than $\frac{w_l}{2}$.

This distance, denoted $d_{(ij)k}$, can be expressed in terms of the seam offsets $x_i$ $x_j$, and $x_k$
\begin{equation}
    \label{eq:lap_stacking_non_parallel}
    \begin{split}
        x_{ij} = \frac{b_i c_j - b_j c_i}{a_i b_j - a_j b_i}, \\
        y_{ij} = \frac{-a_i x_{ij} - c_i}{b_i}, \\
        d_{(ij)k} = |a_k x_{ij} + b_k y_{ij} + c_k|,
    \end{split}
\end{equation}
where $a_*$, $b_*$, $c_*$ are the standard form line parameters associated with seam $*$, with $a_*^2 + b_*^2 = 1$. This distance can be equivalently be expressed as a linear combination of line constants:
\begin{equation}
    \Big{|} \begin{bmatrix}
    \frac{a_j b_k - a_k b_j}{a_i b_j-a_j b_i} & \frac{a_k b_i - a_i b_k}{a_i b_j-a_j b_i} & 1 \\
    \end{bmatrix} 
    \cdot 
    \begin{bmatrix}
    c_i & c_j & c_k
    \end{bmatrix}^T
    \Big{|} 
    \leq d_{(ij)k}^{min},
\end{equation}
where $c_*$ is the line constant of seam $*$ and is a linear function of the seam offset $x_*$.The constraint is linear so long as  the intersection of seams $i$ and $j$ remains in the half-space defined by seam $k$.

The constant minimum distance is:
\begin{equation}
    d_{(ij)k}^{min} = \Big{|}proj_{k_{\perp}} \frac{w_l}{\sin(\theta_{ij}/2)}\Big{|} - \frac{w_l}{2},
\end{equation}
where $\theta_{ij}$ is the interior angle formed by seams $i$ and $j$. For each seam triplet, three constraints are required, one for each combinations $(ij)k$, $(jk)i$, and $(ki)j$.

\begin{figure}
    \centering
    \includegraphics[width=1.0\linewidth]{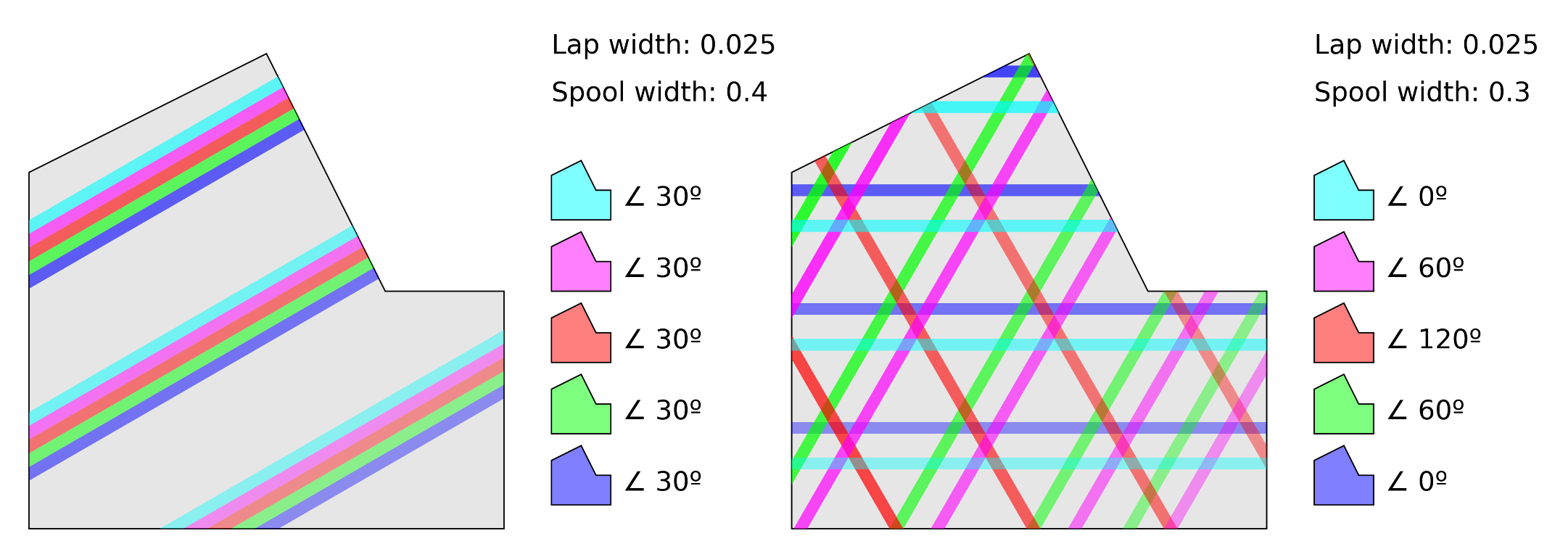}
    \caption{Seam patterns for laminates with single (left) and multiple (right) fiber orientations.}
    \label{fig:example_basic}
\end{figure}

Figure~\ref{fig:example_basic} shows two designs that satisfy the overlap stacking constraint. On the left, with all plies sharing a single fiber orientation, the intersection constraint defined by Eq.~\eqref{eq:lap_stacking_parallel}, together with the ply width constraint in Eq.~\eqref{eq:max_ply_width}, effectively stagger the seams such that overlap stacking is avoided entirely. On the right, with three unique fiber orientations, Eq.~\eqref{eq:lap_stacking_non_parallel} avoids triple overlap stacking, while double overlap stacking remains unavoidable.

Although this formulation guarantees that triple intersections are avoided, it may enforce constraints based on intersections outside the ply polygon, thereby over constraining the design space. To circumvent this issue, intersections outside the ply polygon are ignored. If a design update results in a related triple overlap inside the design polygon, the set of constraints for the next seam placement will rectify the violation. This solution requires that the local optimization be rerun on the final solution to ensure that it is indeed valid.

\subsection{Stay-out zones}
\begin{figure}[h]
    \centering
    \includegraphics[width=1.0\linewidth]{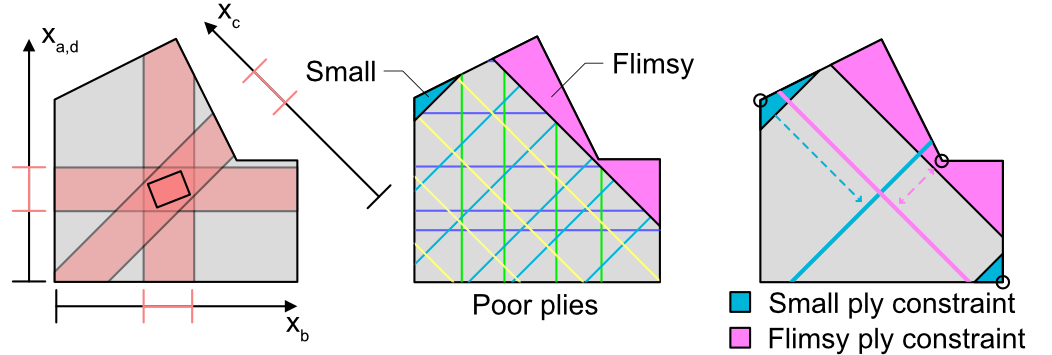}
    \caption{Left: A stay-out zone and its projection onto the design space corresponding to ply orientations defined in Fig.~\ref{fig:overlap_illustration}. Center: Example seam pattern with poor quality sub-plies highlighted. Right: Fiber orientation-dependent minimum distance constraints (dashed lines) targeting small and flimsy ply conditions.}
    \label{fig:additional_constraints}
\end{figure}

Stay-out zones segment the design space into allowed and prohibited domains. Consider the example in Fig.~\ref{fig:additional_constraints}. For each ply, the projection of the stay-out polygon $S$ onto a coordinate $d_{\perp}^j$ through the ply origin and normal to the fiber direction of ply $j$ defines a domain of prohibited offset values. Mathematically,
\begin{equation}
    x_*^j \notin proj_{d_{\perp}^j} S,
\end{equation}
where $x_*^j$ is the offset of any seam on ply $j$.

\begin{figure}[h]
    \centering
    \includegraphics[width=1.0\linewidth]{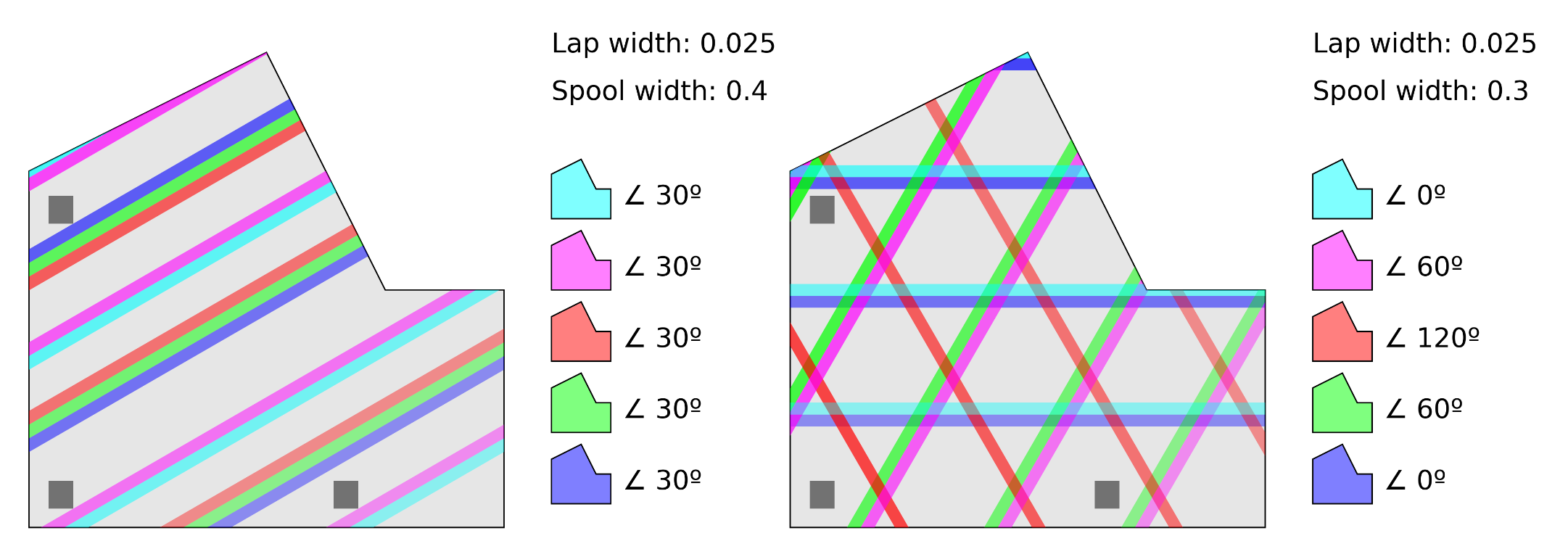}
    \caption{Seam patterns for laminates with single (left) and multiple (right) fiber orientations with stay-out zones (dark gray).}
    \label{fig:example_stayouts}
\end{figure}

Figure~\ref{fig:example_stayouts} shows how the added stay-out zones (gray) would affect the splicing of the laminates from Fig.~\ref{fig:example_basic}.

\subsection{Ply quality}
An additional consideration with respect to seam placement relates to the size and shape of the resulting sub-plies. Specifically, narrow, small, or flimsy sub-plies such as those shown in Fig.~\ref{fig:additional_constraints} should be avoided to ensure that they do not drift from their intended position during fabrication.

These requirements can be expressed as distance constraints to specific boundary vertices, as shown in Fig.~\ref{fig:additional_constraints} (right). The pertinent vertices for avoiding small sub-plies $V_{small}^j$ on ply $j$ are those whose edges form a cone pointing towards the seam. Formally:
\begin{equation}
    \begin{split}
        V_{small}^j =  & \{{} v | v \in V, \angle e_{i} e_{i+1} \in [0, \ \pi]  \\
        & \wedge \sign(proj_{d_\perp} \ e_i) = \sign(proj_{d_\perp^j} \ e_{i+1}) \}
    \end{split}
\end{equation}
where $V$ is the set of all $n$ ply polygon vertices, $e_i$ and $e_{i+1}$ are neighboring polygon edges, and $\mathbf{d}_\perp^j$ is the fiber-transverse direction.
Similarly, the pertinent vertices for flimsy sub-plies $V_{flim}^j$ are those which form a cone pointing away from the seam:
\begin{equation}
    \begin{split}
        V_{flim}^j = & \{ v | v \in V, \angle e_{i} e_{i+1} \in (\pi, \ 2\pi) \\
        & \wedge \sign(proj_{d_\perp} \ e_i) = \sign(proj_{d_\perp^j} \ e_{i+1}) \} 
        \end{split}
\end{equation}

In each case, an inequality constraint can be added with a user-provided minimum value:
\begin{equation}
    \begin{split}
        |x_*^j - x(V_{k \ flim}^j)| \geq x_{flim}^{min} \\
        |x_*^j - x(V_{k \ small}^j)| \geq x_{small}^{min},
    \end{split}
\end{equation}
where $x(V_{k \ flim}^j)$ is the offset of ply polygon vertex $k$ in $V_{flim}^j$, and $x(V_{k \ small}^j)$ is the offset of ply polygon vertex $k$ in $V_{small}^j$, respectively.

\begin{figure}[h]
    \centering
    \includegraphics[width=1.0\linewidth]{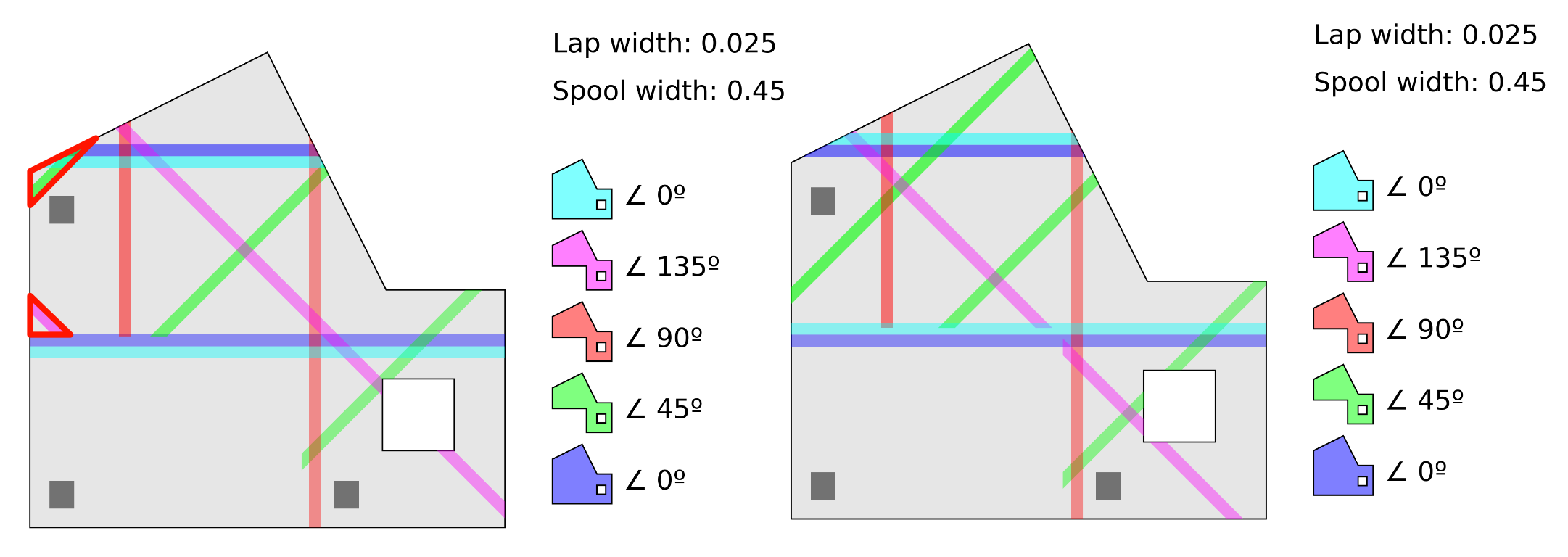}
    \caption{Seam patterns generated without (left) and with (right) ply quality constraints. Poor-quality sub-plies are contoured in red.}
    \label{fig:example_ply_quality}
\end{figure}

Figure~\ref{fig:example_ply_quality} shows how the inclusion of ply quality constraints alters the seam pattern. Without the ply quality constraints, the green and magenta sub-plies contain small corner sub-plies. With the constraints activated, all sub-plies have a minimum width greater than the user-set value of $w_s/2$.

\subsection{Optimization}
\label{sec:optimization}
With the manufacturing requirements formulated as constraints, we devise an optimization problem for which the objective is to partition each ply into producible sub-plies, using as few seams as possible. If seams are placed one at a time, the local objective is to maximize the coverage for a set of existing seams:
\begin{equation}
\label{eq:strip_width_max}
\begin{split}
    F & = -\sum_{j=1}^m \sum_{i=1}^{n^j-1} (x_{i+1}^j - x_i^j) \\
    & = -\sum_{j=1}^m x_{n}^j
\end{split}
\end{equation}
where $m$ is the number of plies, $x_i^j$ is the design variable associated with the $i$\textsuperscript{th} seam on ply $j$, and $n^j$ is the number of seams on ply $j$. The $x_1^j$ term drops out since it is always set to align with the ply origin and is therefore equal to 0.

With the constraints described previously, the optimization problem is as follows,
\begin{equation}
\begin{aligned}
\minimize \limits_{\mathbf{x}} \quad & -\sum_{j=1}^m x_{n}^j&\\
\textrm{s.t.} \quad & w_{min} \leq x_{i+1}^j - x_{i}^j \leq w_t - w_l, & \forall \ i, j \\
& |x_j - x_i| \geq w_l & \forall s_i^* \parallel s_j^* \\
& d_{(ij)k} \leq d_{(ij)k}^{min} & \forall \ s_i^* \nparallel s_j^* \nparallel s_k^*\\
& x_*^j \notin proj_{d_{\perp}^j} S & \forall j, S \\
& |x_*^j - x(V_{k \ flim}^j)| \geq x_{flim}^{min} & \forall j, k \\
& |x_*^j - x(V_{k \ small}^j)| \geq x_{small}^{min} & \forall j, k \\
& x_i^j \geq 0, & \forall \ i, j \\
\end{aligned}
\label{eq:opt}
\end{equation}

\begin{figure}
    \centering
    \includegraphics[width=1.0\linewidth]{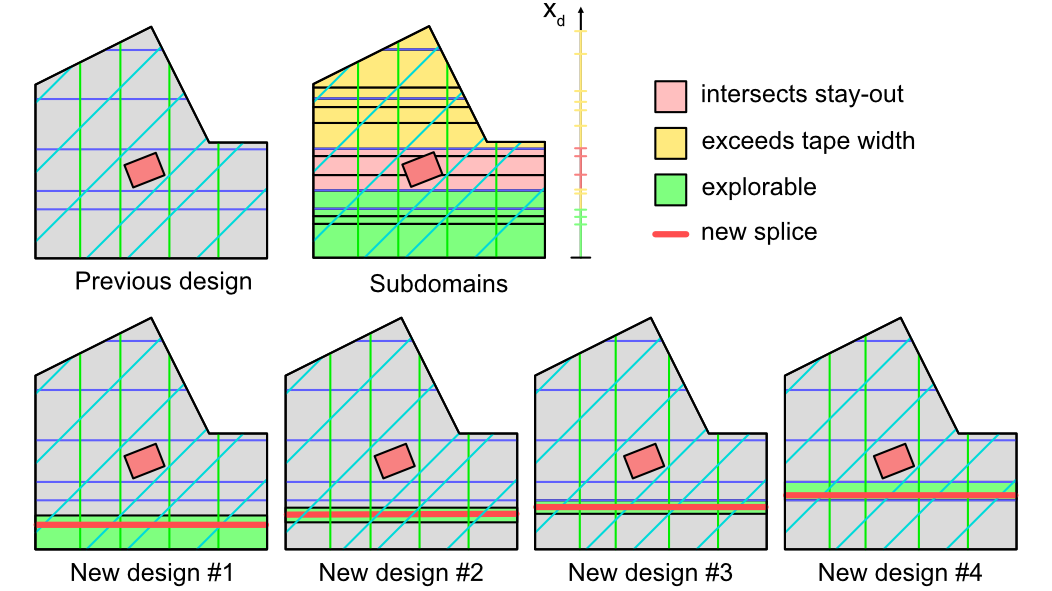}
    \caption{Seam insertion procedure illustrating how the design domain is first divided into linear subdomains based on the current design, then each is locally optimized.}
    \label{fig:local_optimization_illustration}
\end{figure}

\begin{figure}
    \centering
    \includegraphics[width=1.0\linewidth]{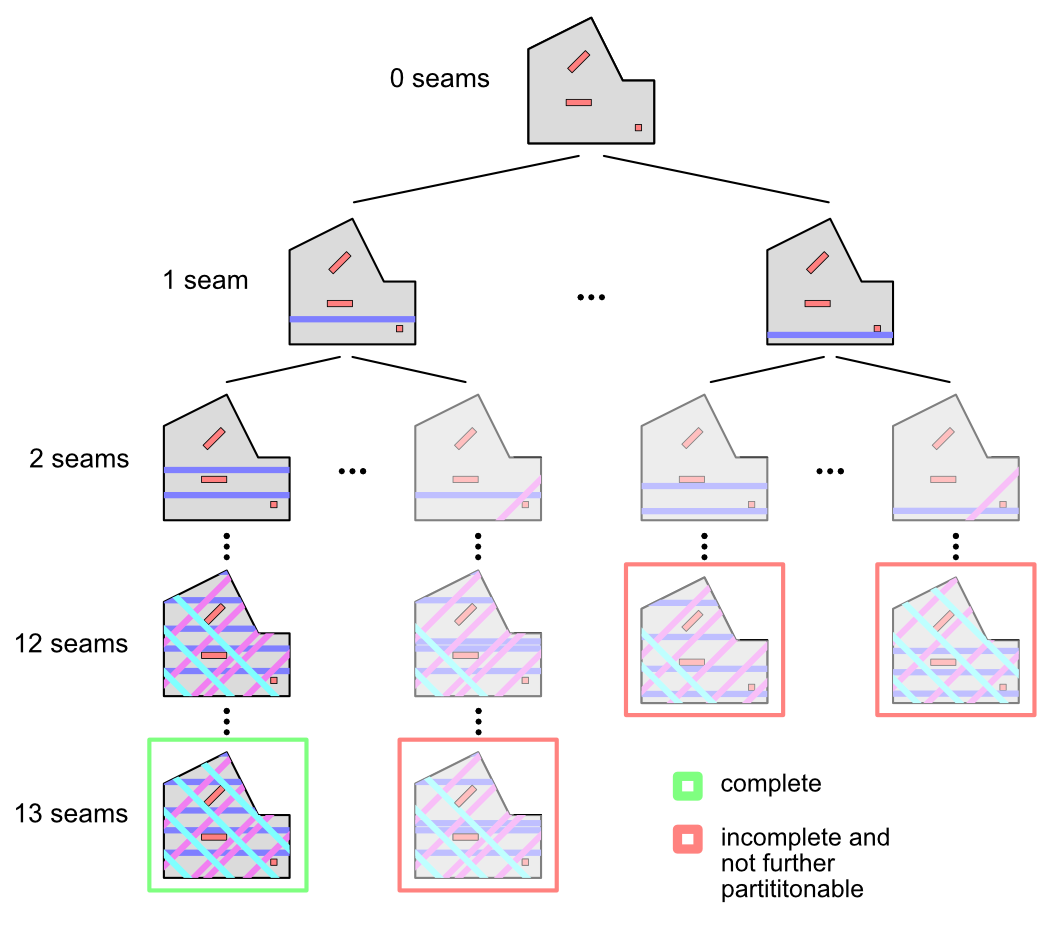}
    \caption{Illustration of the beam search strategy. Each row of the design tree represents all designs with one additional seam from the previous.}
    \label{fig:beam_search_illustration}
\end{figure}

For a given set of seams, the optimization problem, as formulated, is piecewise linear. Therefore, this local problem can be solved efficiently via linear programming. The global problem in which any viable design must also fully partition all plies is discrete and, at first glance, cannot be solved as efficiently. For this, one could adopt an inelegant but robust beam search strategy. Specifically, a design tree in which each node represents a locally optimized set of seams is explored breadth first placing seams in each linear-subspace, keeping a fixed number of best performing designs from each generation. 
The procedure terminates when one of the designs fully partitions all plies. The search strategy is presented graphically in Figs.~\ref{fig:local_optimization_illustration} and \ref{fig:beam_search_illustration}. 

In the parameter study presented in Section~\ref{sec:efficacy}, it appears that a greedy search always arrives at the same final design as the beam search, regardless of the beam width. Thus, it appears that the beam search is unnecessary and the efficient greedy strategy is sufficient. The greedy solution is presented algorithmically in Alg.~\ref{alg:design}.

\RestyleAlgo{ruled}
\begin{algorithm}[t]
\caption{Generate seam pattern}\label{alg:design}
\KwData{$\{P\}, \{S\}, w_s, w_l, w_{min}$}
\KwResult{Set of seams $\mathbf{x}$}
$\mathbf{x} \gets \{\}$ \\
\For{each ply in \{P\}}{
    \While{last seam inside ply polygon}{
        Identify linear sub-spaces \\
        \While {seam not placed}{
            \For{each linear sub-space, starting from farthest to previous seam}{
                Initialize new seam in sub-space \\
                Optimize all existing seams \\
                \If {optimization success}{
                    $\mathbf{x} \gets solution$
                }
            }
        }
    }
}
\end{algorithm}

\subsection{Cost minimization}
\begin{figure}
    \centering
    \includegraphics[width=1.0\linewidth]{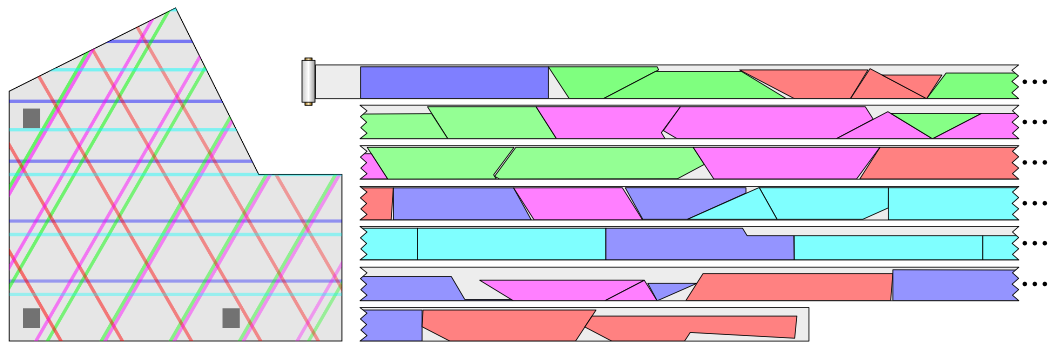}
    \caption{Example seam pattern and corresponding prepreg cutting pattern (presented in pieces). Trim loss (gray) is concentrated above narrow sub-plies.}
    \label{fig:nest}
\end{figure}

In addition to generating viable designs given a prescribed spool width, engineers may wish to explore multiple spool width options and evaluate the trade offs in terms of design fidelity and manufacturing costs. This section describes an alternative formulation in which the objective function is an estimation of total manufacturing costs. The following linear cost function is proposed:
\begin{equation}
    \label{eq:cost}
    C :=  A_{mat} \ \hat{c}_{mat} + n_{seam} \ \hat{c}_{seam},
\end{equation}
where the first term relates to material cost, and the second one relates to processing cost. The former accounts for the total cost of consumed material, including overlaps and material wastage, while the latter accounts for labor and consumables, which are largely proportional to the number of seams. 
We assume that material cost per unit area $\hat{c}_{mat}$ is a linear function of the spool width. Further, we divide the total material area $A_{mat}$ into design area $A_d$, overlap area $A_l$, and trim loss $A_{trim}$, which is material wasted when ply patterns are cut from the prepreg spool.
$A_d$ is a constant value, defined by the original ply layup design. $A_l$ is a function of the number of seams, the overlap width, and the position of each seam within the part geometry. Here, since the polygon defined by each overlap has length much greater than width, its area can be approximated as that of a rectangle with width equal to the user-defined overlap width. Furthermore, since seams are spread relatively evenly over each ply, the length component can be approximated as constant for each ply. $A_l$ can therefore be approximated as a linear function of the number of seams. Finally, the trim loss term $A_{trim}$ is an estimation of the leftover material after sub-plies have been nested and cut from the spool stock. Since the cutting pattern is generally not a continuous function of the seam offsets, a simple approximation is proposed based on the observation that trim loss tends to be positively correlated with the difference between the width of each sub-ply and the spool width. Figure~\ref{fig:nest} shows an example nesting pattern generated by placing sub-plies sequentially, each time selecting the remaining sub-ply that leaves the smallest gap with the previous one. By considering the 180$^{\circ}$ rotated and mirrored configurations of each sub-ply, this simple strategy produces reasonably compact nesting patterns, without iteration or heavy computation. Notice that the trim loss is concentrated in two areas: (1) at interfaces between horizontally adjacent sub-plies, and (2) above narrower sub-plies that do not extend to the edge of the spool. In real-world applications, where the sub-ply length-to-width ratio is typically greater, the material waste above narrow sub-plies is expected to account for the majority of the trim loss.

If the sub-plies are approximated as rectangles of length $\Bar{h}_sp$, then the trim loss term can be written as 
\begin{equation}
     A_{trim} = \sum_{j=1}^m \sum_{i=1}^{n^j-1} w_s - (x_{i+1}^j - x_i^j) \ \Bar{h}_{sp}
\end{equation}

Finally, the discrete $n_{seam}$ term in Eq.~\eqref{eq:cost} can be replaced by the continuous function defined in Eq.~\eqref{eq:strip_width_max}.
Eq.~\eqref{eq:cost} is then fully linearized as: 

\begin{equation}
    \label{eq:reduced_cost}
    \begin{split}
         C = (A_d + \sum_{j=1}^m \sum_{i=1}^{n^j-1} w_s - (x_{i+1}^j - x_i^j) \ \Bar{h}_{sp}) \cdot (\hat{a}_{mat} \cdot w_s + \hat{b}_{mat}) \\
         + \hat{c}_{seam} \ \sum_{j=1}^m w_s \cdot (n^j-1) - x_{n}^j
    \end{split}
\end{equation}
where $\hat{a}_{mat}$ and $\hat{b}_{mat}$ are the linear factor and constant relating the material unit cost to the spool width. The solution strategy is nearly identical to the original, with the addition of user-defined minimum $w_s^{min}$ and maximum spool width $w_s^{max}$ parameters, and substitution of the linear solver with a quadratic solver~\cite{nocedal2006quadratic}.

\section{Results}
\label{sec_results}

The proposed partitioning strategy was implemented in python, using the Numpy module, and the least squares SLSQP method from the Scipy Optimize module. Results were produced on a MacBook Pro with M1 Max and 64 GB integrated memory. Compute time for each experiment was on the order of seconds.

\subsection{Wing example}
As a representative example, the partitioning strategy was applied to the design of an airplane wing, and to the body panels of an armored vehicle. In addition to producing viable designs for a given set of user inputs, the results presented herein illustrate how designers can use the partitioning algorithm to quickly explore the implications of each parameter on material utilization, cost, and design fidelity. 

In these examples, an exaggerated overlap width of 1.2\% (normalized with respect to the flat part geometry) was used, for clarity. In industry, nominal overlap widths typically vary between 2 and 10 cm, and are dictated by engineering specification and process parameters such as prepreg thickness, part geometry, and tooling design. Note that the exaggerated overlap width restricts the solution space, and reduces the quality of the results.

\begin{figure*}[h]
    \centering
    \includegraphics[width=1.0\linewidth]{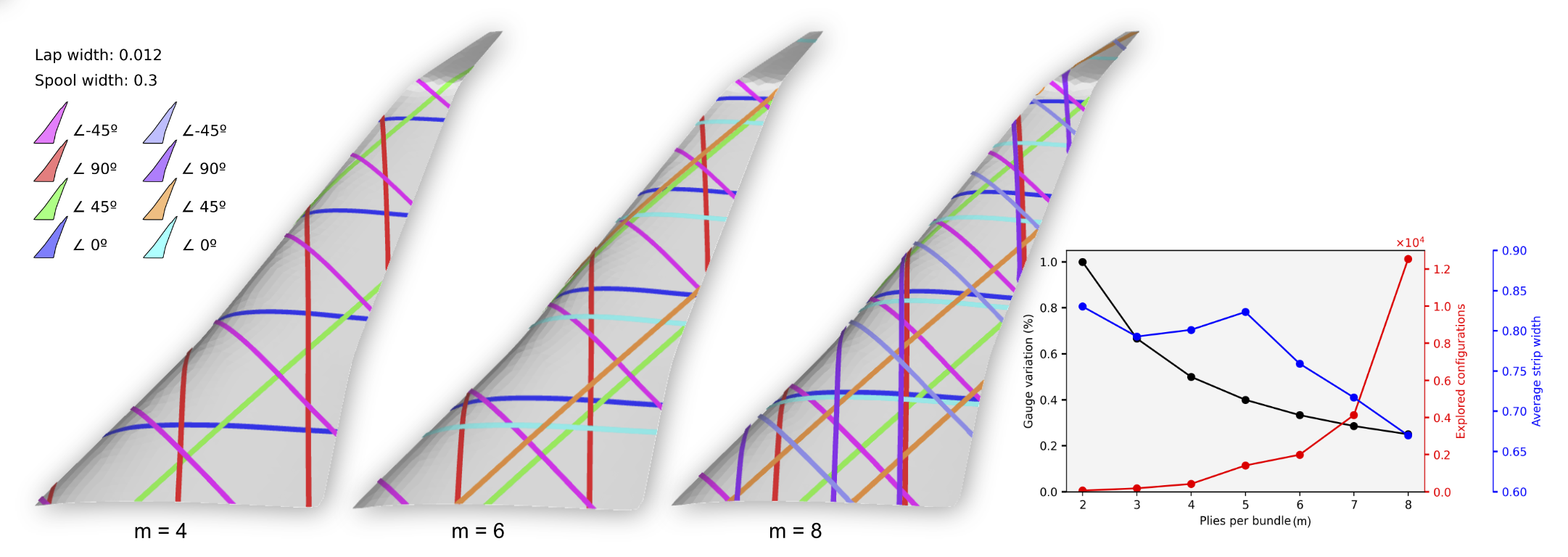}
    \caption{Airplane wing and seam patterns obtained for bundles of $m =$ 4, 6, and 8 plies, with fiber orientations as shown. The prepreg spool width $w_s$ and overlap width $w_{l}$ were 0.3 and 1.2\%, respectively, normalized with respect to the wing length. Performance in terms of laminate thickness variation, average ply width, and the number of explored nodes are plotted on the right.}
    \label{fig:wing_example_1}
\end{figure*}

In Fig.~\ref{fig:wing_example_1}, the algorithm is applied to a generic airplane wing with thickness tolerance between 100\% (left) and 25\% (right). Internally, this is achieved by considering between 2 and 8 plies per bundle. 

As the thickness tolerance is tightened, the problem is increasingly constrained, and it is more difficult to ensure that each ply width is maximized. Predictably, the resulting partitioning patterns include more seams and smaller sub-plies as the number of plies per bundle $m$ is increased. 

\subsection{Time complexity and ply sorting}
The local linear programming problem is solved efficiently using the standard simplex method~\cite{ficken2015simplex}. Though it is not strictly polynomial, other linear programming methods are known to be~\cite{megiddo1986complexity}.

The greedy search has time complexity proportional to the number of explored nodes, which itself depends on the time to formulate and the resolvability of the local problem. In practice, this is a function of the number of seams already placed and the number of seam intersections.
To improve performance, the latter can be reduced by adjusting the ordering of plies. Specifically, since large component layups are typically made up of many similarly oriented plies, they can be sorted by orientation, and processed as bundles containing a minimum of unique fiber orientations. This reduces the number of seam intersections, and therefore also reduces the number of linear sub-spaces, and improves local problem resolvability. If possible, only a single orientation should be included in each bundle. Otherwise, selecting two orientations that are as perpendicular as possible is preferable since it minimizes the total overlap area.

\begin{figure*}[t]
    \centering
     \includegraphics[width=1.0\linewidth]{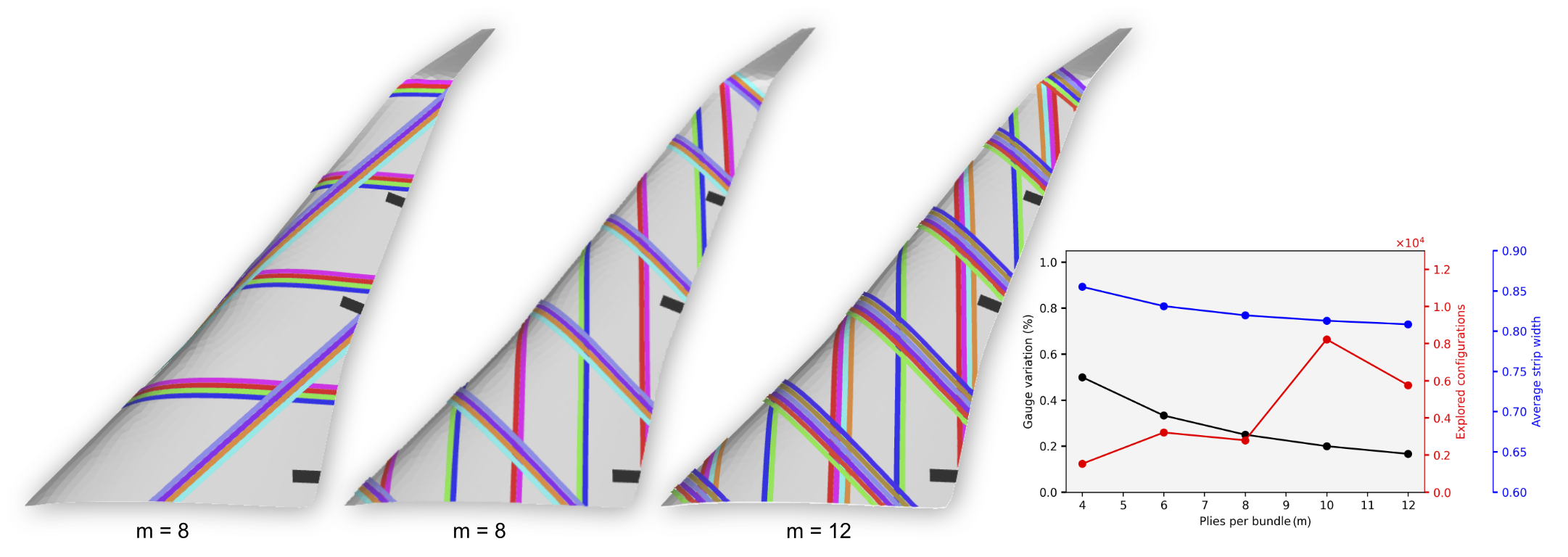}
    \caption{Seam patterns obtained for bundles of similarly-oriented plies. The spool width and overlap width were 0.3 and 1.2\%, normalized with respect to the wing length. Stay-outs are shown in black.}
    \label{fig:wing_example_2}
\end{figure*}

To illustrate this optimization, consider that the seam patterns in Fig.~\ref{fig:wing_example_1} correspond to the first 8 plies of a 24-ply layup. If the layup consists of 6 plies for each orientation, better performance can be achieved by processing ply bundles made up of only one or two unique orientations. In Fig.~\ref{fig:wing_example_2}, the seam patterns resulting from 3 bundles of 8 plies, and 2 bundles of 12 plies are shown. In comparison to the 8-ply bundle from Fig.~\ref{fig:wing_example_1}, the number of nodes traversed is considerably lower, even for the larger bundles. Moreover, the resulting seam patterns are of higher quality since all sub-plies have maximal width (excluding end plies). This is true even for the larger 12-ply bundles which achieve an improved thickness variation tolerance of 17\%. 

\subsection{Global search efficacy}
\label{sec:efficacy}

\begin{figure}[t]
    \centering
     \includegraphics[width=1.0\linewidth]{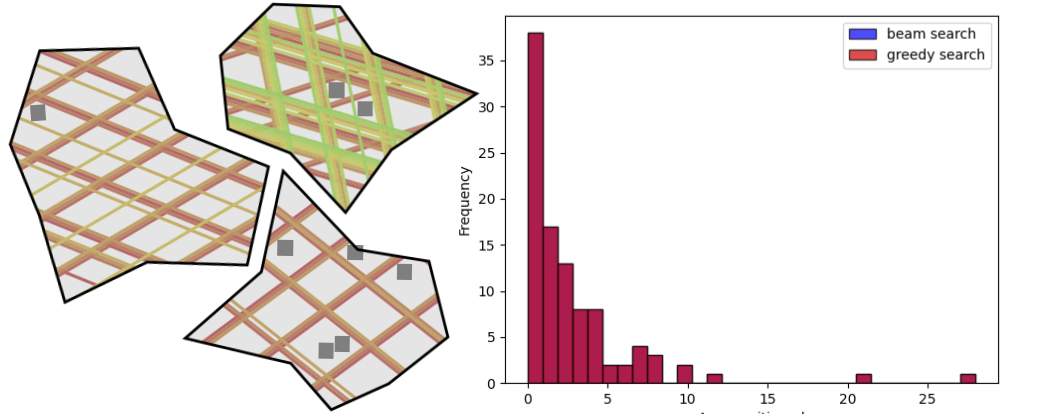}
    \caption{Left: example synthetic problems and final seam patterns. Right: overlaid histograms of total ply area partitioned prior to failure via greedy search and beam search.}
    \label{fig:beam_width}
\end{figure}

To ascertain whether the greedy search is sufficient for the global optimization problem, we compare the results to those obtained via beam search. The solver was run until failure (not further partitionable) using synthetic geometry consisting of an infinite number of randomly generated plies and 1-20 randomly positioned square stay-out zones. Solver parameters were $w_s = 0.2$, $w_l = 0.01$, $w_{min}=0.1$. Examples are shown in Fig.~\ref{fig:beam_width}, along with a histogram of the total partitioned length. 
In the beam search, a beam width of 10~000 was used to guarantee the exploration of all design nodes, effectively making the search exhaustive. Remarkably, running the solver 100~000 times revealed no performance improvement whatsoever, and no difference in the final designs. In all cases, this corresponded to an exhaustive search.

We note, however, that the beam search may still be advantageous for more complex design objectives, such as those including structural analysis or other non-convex terms.

\subsection{Vehicle example}
\begin{figure*}[h]
    \centering
    \includegraphics[width=1.0\linewidth]{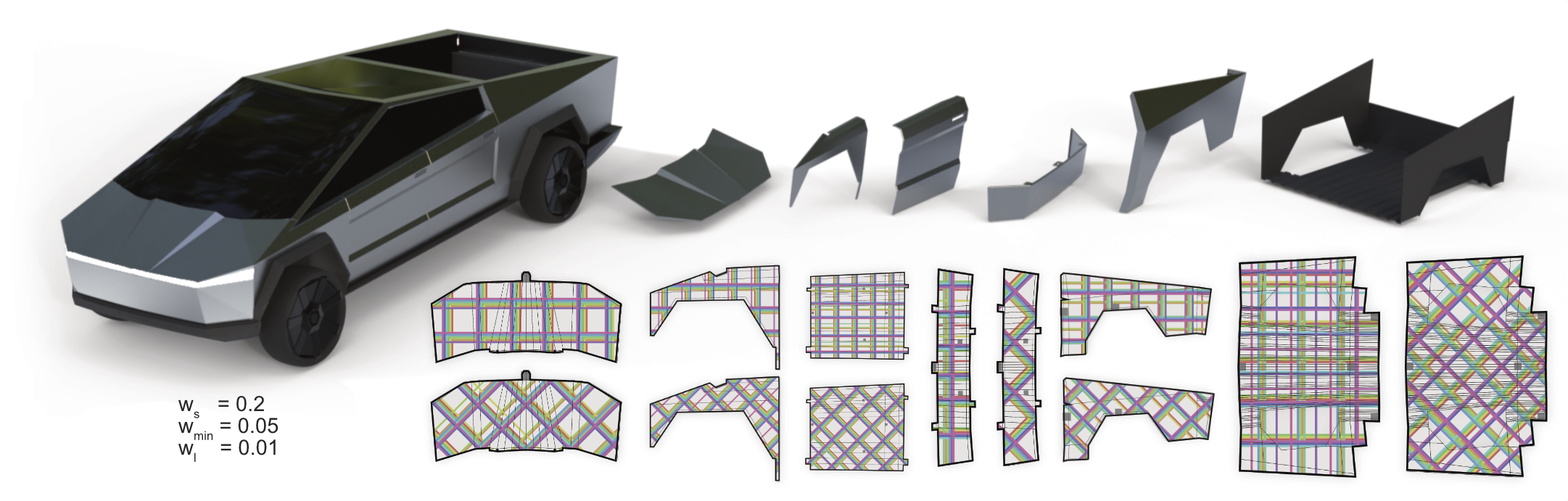}
    \caption{Retro futuristic armored vehicle body panels and their respective seam patters for two-orientation ply bundles.}
    \label{fig:truck_example}
\end{figure*}

As a second example, to illustrate its practical utility, the proposed method is applied to the body panels of a retro-futuristic armored vehicle. Each of the panels can be fabricated from UD prepreg, using the partitioning method to minimize defects. In Fig.~\ref{fig:truck_example}, 24-ply layups were partitioned with prescribed stay-out zones, spool width of 20\%, normalized to the largest dimension of the axis-aligned geometry bounding box, lap width of 1\% and minimum ply width 5\%. Perpendicular fiber orientations were bundled for efficiency, achieving worst-case tolerance of 16.7\%. For each panel, viable seam patterns were identified within 20~seconds. Note that the exaggerated overlap width was the limiting factor, and much larger ply bundles can be handled with more conservative overlap width.

\subsection{Design exploration}
\begin{figure*}[t]
    \centering
    \includegraphics[width=1.0\linewidth]{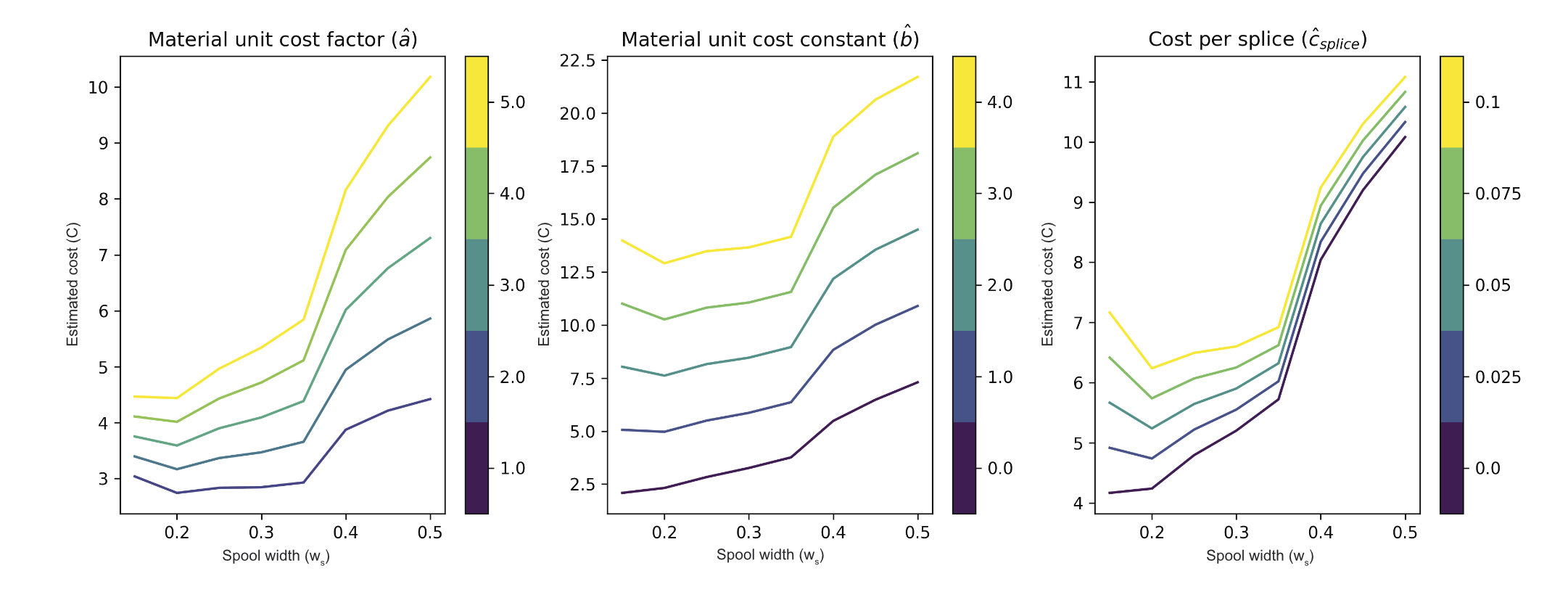}
    \caption{Estimated production cost as a function of spool width. In each plot, $\hat{a} = 5.0$, $\hat{b} = 1.0$, and $\hat{c}_{seam} = 0.01$ unless otherwise specified.}
    \label{fig:wing_example_3}
\end{figure*}

In the following example, we demonstrate how the proposed method can be used to explore design trade-offs by plotting the estimated production cost with respect to the manufacturing constraints. Using the same wing model, the total production cost is estimated for a $[0^{\circ}, 60^{\circ}, 120^{\circ}]$ ply bundle over a range of spool widths. Figure~\ref{fig:wing_example_3} illustrates the impact of the material unit cost parameter, as well as the manufacturing cost per seam on the overall cost. 

Naturally, as the spool width decreases, the number of seams in the resulting seam pattern increases, as does material usage and final part weight. However, since spool width is positively correlated with material unit cost, the overall cost can be minimized by finding the optimal balance between spool width and material usage. 

\section{Concluding remarks}
\label{sec_conclusions}

The method presented in this work offers a robust and efficient alternative to manual ply partitioning for large-scale laminar composite manufacturing. By leveraging the developable nature of the detailed ply layups produced in the digital process planning pipeline, the partitioning problem is linearized and solved through a combination of linear programming and heuristic tree search. The proposed algorithm allows users to specify various manufacturing constraints, as well as to explore design options and trade-offs in a matter of minutes, without the need for tedious and unreliable trial-and-error. 

\subsection{Limitations}
While the proposed method can be applied to a wide range of design applications, it is best suited to convex geometry, and may be overly conservative for non-convex problems in which seams intersect multiple boundaries. This could be addressed by subdividing the ply polygon into convex regions. In addition, this method is limited to straight fiber-aligned seams and developable surfaces. The inclusion of fiber-transverse and staircase seams would expand the design space and should be considered in future work. However, extension to non-developable geometry is likely not possible.
Future work should also address special cases in which layup designs cannot be straightforwardly discretized into flattenable zones, but which nonetheless can be constructed from flat prepreg plies. The proposed method could be extended for such cases. 
Finally, the proposed method does not directly consider the mechanical performance of the as-manufactured part. In other words, we are not explicitly modeling the effects of manufacturing-induced defects. Our approach could potentially be coupled with structural analysis to better estimate design fidelity and structural integrity. In such case, the greedy search may not be sufficient for the global optimization, and the beam search may offer better performance.
%-------------------------------------------------------------------------

% bibtex
\bibliographystyle{eg-alpha-doi} 
\bibliography{references}       

\newcommand{\etalchar}[1]{$^{#1}$}
\begin{thebibliography}{\uppercase{SFAKP02}}

\bibitem[AAD98]{alexander1998part}
\textsc{Alexander P., Allen S., Dutta D.}:
\newblock Part orientation and build cost determination in layered manufacturing.
\newblock \emph{Computer-Aided Design 30}, 5 (1998), 343--356.

\bibitem[AHL17]{alexa2017optimal}
\textsc{Alexa M., Hildebrand K., Lefebvre S.}:
\newblock Optimal discrete slicing.
\newblock \emph{ACM Transactions on Graphics (TOG) 36}, 1 (2017), 1--16.

\bibitem[BL09]{bruno2009you}
\textsc{Bruno V., L{\'e}vy B.}:
\newblock What you seam is what you get.

\bibitem[CDMS06]{canellidis2006pre}
\textsc{Canellidis V., Dedoussis V., Mantzouratos N., Sofianopoulou S.}:
\newblock Pre-processing methodology for optimizing stereolithography apparatus build performance.
\newblock \emph{Computers in industry 57}, 5 (2006), 424--436.

\bibitem[CLRS22]{cormen2022introduction}
\textsc{Cormen T.~H., Leiserson C.~E., Rivest R.~L., Stein C.}:
\newblock \emph{Introduction to algorithms}.
\newblock MIT press, 2022.

\bibitem[DTS16]{delfs2016optimized}
\textsc{Delfs P., Tows M., Schmid H.-J.}:
\newblock Optimized build orientation of additive manufactured parts for improved surface quality and build time.
\newblock \emph{Additive Manufacturing 12} (2016), 314--320.

\bibitem[EC93]{elber1993tool}
\textsc{Elber G., Cohen E.}:
\newblock Tool path generation for freeform surface models.
\newblock In \emph{Proceedings on the second ACM symposium on Solid modeling and applications} (1993), pp.~419--428.

\bibitem[EF97]{elber19975}
\textsc{Elber G., Fish R.}:
\newblock 5-axis freeform surface milling using piecewise ruled surface approximation.

\bibitem[Fic15]{ficken2015simplex}
\textsc{Ficken F.~A.}:
\newblock \emph{The simplex method of linear programming}.
\newblock Courier Dover Publications, 2015.

\bibitem[Gra24]{Grand_View_Research_2024}
2024.
\newblock URL: \url{https://www.grandviewresearch.com/industry-analysis/laminar-composites-market-report}.

\bibitem[JDM{\etalchar{*}}17]{jin2017optimization}
\textsc{Jin Y., Du J., Ma Z., Liu A., He Y.}:
\newblock An optimization approach for path planning of high-quality and uniform additive manufacturing.
\newblock \emph{The International Journal of Advanced Manufacturing Technology 92} (2017), 651--662.

\bibitem[KEBP15]{kim2015precise}
\textsc{Kim Y.-J., Elber G., Barto{\v{n}} M., Pottmann H.}:
\newblock Precise gouging-free tool orientations for 5-axis cnc machining.
\newblock \emph{Computer-Aided Design 58} (2015), 220--229.

\bibitem[KZW{\etalchar{*}}15]{kwok2015styling}
\textsc{Kwok T.-H., Zhang Y.-Q., Wang C.~C., Liu Y.-J., Tang K.}:
\newblock Styling evolution for tight-fitting garments.
\newblock \emph{IEEE transactions on visualization and computer graphics 22}, 5 (2015), 1580--1591.

\bibitem[LBRM12]{luo2012chopper}
\textsc{Luo L., Baran I., Rusinkiewicz S., Matusik W.}:
\newblock Chopper: Partitioning models into 3d-printable parts.
\newblock \emph{ACM Transactions on Graphics (TOG) 31}, 6 (2012), 1--9.

\bibitem[M{\etalchar{*}}86]{megiddo1986complexity}
\textsc{Megiddo N., et~al.}:
\newblock \emph{On the complexity of linear programming}.
\newblock Citeseer, 1986.

\bibitem[MGVL19]{mehdikhani2019voids}
\textsc{Mehdikhani M., Gorbatikh L., Verpoest I., Lomov S.~V.}:
\newblock Voids in fiber-reinforced polymer composites: A review on their formation, characteristics, and effects on mechanical performance.
\newblock \emph{Journal of Composite Materials 53}, 12 (2019), 1579--1669.

\bibitem[MRI00]{masood2000part}
\textsc{Masood S.~H., Rattanawong W., Iovenitti P.}:
\newblock Part build orientations based on volumetric error in fused deposition modelling.
\newblock \emph{The International Journal of Advanced Manufacturing Technology 16}, 3 (2000), 162--168.

\bibitem[MTMP20]{montes2020computational}
\textsc{Montes J., Thomaszewski B., Mudur S., Popa T.}:
\newblock Computational design of skintight clothing.
\newblock \emph{ACM Transactions on Graphics (TOG) 39}, 4 (2020), 105--1.

\bibitem[MWJ12]{meng2012flexible}
\textsc{Meng Y., Wang C.~C., Jin X.}:
\newblock Flexible shape control for automatic resizing of apparel products.
\newblock \emph{Computer-aided design 44}, 1 (2012), 68--76.

\bibitem[noa23]{noauthor_catia_2023}
{CATIA} {V5}, May 2023.
\newblock URL: \url{https://www.3ds.com/products/catia/catia-v5}.

\bibitem[NSLV86]{nee1986designing}
\textsc{Nee A., Seow K., Long S., Venkatesh V.}:
\newblock Designing algorithm for nesting irregular shapes with and without boundary constraints.
\newblock \emph{CIRP Annals 35}, 1 (1986), 107--110.

\bibitem[NW06]{nocedal2006quadratic}
\textsc{Nocedal J., Wright S.~J.}:
\newblock Quadratic programming.
\newblock \emph{Numerical optimization} (2006), 448--492.

\bibitem[PDF{\etalchar{*}}22]{pietroni2022computational}
\textsc{Pietroni N., Dumery C., Falque R., Liu M., Vidal-Calleja T.~A., Sorkine-Hornung O.}:
\newblock Computational pattern making from 3d garment models.
\newblock \emph{ACM Trans. Graph. 41}, 4 (2022), 157--1.

\bibitem[PTH{\etalchar{*}}17]{poranne2017autocuts}
\textsc{Poranne R., Tarini M., Huber S., Panozzo D., Sorkine-Hornung O.}:
\newblock Autocuts: simultaneous distortion and cut optimization for uv mapping.
\newblock \emph{ACM Transactions on Graphics (TOG) 36}, 6 (2017), 1--11.

\bibitem[PW99]{pottmann1999approximation}
\textsc{Pottmann H., Wallner J.}:
\newblock Approximation algorithms for developable surfaces.
\newblock \emph{Computer Aided Geometric Design 16}, 6 (1999), 539--556.

\bibitem[RCF{\etalchar{*}}19]{ren2019thermo}
\textsc{Ren K., Chew Y., Fuh J., Zhang Y., Bi G.}:
\newblock Thermo-mechanical analyses for optimized path planning in laser aided additive manufacturing processes.
\newblock \emph{Materials \& Design 162} (2019), 80--93.

\bibitem[SE06]{subag2006piecewise}
\textsc{Subag J., Elber G.}:
\newblock Piecewise developable surface approximation of general nurbs surfaces, with global error bounds.
\newblock In \emph{International conference on geometric modeling and processing} (2006), Springer, pp.~143--156.

\bibitem[SFAKP02]{smith2002optimal}
\textsc{Smith T.~S., Farouki R.~T., Al-Kandari M., Pottmann H.}:
\newblock Optimal slicing of free-form surfaces.
\newblock \emph{Computer Aided Geometric Design 19}, 1 (2002), 43--64.

\bibitem[Sof24]{fibersim_siemens}
\textsc{Software S. D.~I.}:
\newblock Fibersim for composite design and manufacturing.
\newblock \url{https://plm.sw.siemens.com/en-US/nx/products/fibersim-composites/}, 2024.

\bibitem[Sta15]{staab2015laminar}
\textsc{Staab G.}:
\newblock \emph{Laminar composites}.
\newblock Butterworth-Heinemann, 2015.

\bibitem[TM22]{tura2022characterization}
\textsc{Tura A.~D., Mamo H.~B.}:
\newblock Characterization and parametric optimization of additive manufacturing process for enhancing mechanical properties.
\newblock \emph{Heliyon 8}, 7 (2022).

\bibitem[TPR04]{thrimurthulu2004optimum}
\textsc{Thrimurthulu K., Pandey P.~M., Reddy N.~V.}:
\newblock Optimum part deposition orientation in fused deposition modeling.
\newblock \emph{International Journal of Machine Tools and Manufacture 44}, 6 (2004), 585--594.

\bibitem[WHZ{\etalchar{*}}23]{wolff2023designing}
\textsc{Wolff K., Herholz P., Ziegler V., Link F., Br{\"u}gel N., Sorkine-Hornung O.}:
\newblock Designing personalized garments with body movement.
\newblock In \emph{Computer Graphics Forum} (2023), vol.~42, Wiley Online Library, pp.~180--194.

\bibitem[WSY02]{wang2002surface}
\textsc{Wang C.~C., Smith S.~S., Yuen M.~M.}:
\newblock Surface flattening based on energy model.
\newblock \emph{Computer-Aided Design 34}, 11 (2002), 823--833.

\end{thebibliography}
%-------------------------------------------------------------------------
\end{document}